\documentclass[12pt]{article}
\usepackage{latexsym}
\usepackage{amsmath}
\usepackage{amssymb}
\newtheorem{thm}{Theorem}[section]

\newtheorem{prop}[thm]{Proposition}
\newtheorem{exam}[thm]{Example}
\newtheorem{cor}[thm]{Corollary}
\newcommand{\Cnt}{\mathrm{Cnt}\,}
\newcommand{\Cntr}{\mathrm{Cnt_r}\,}
\newcommand{\Int}{\mathrm{Int}\,}
\newcommand{\Aut}{\mathrm{Aut}\,}
\newcommand{\Ad}{\mathrm{Ad}\,}
\newcommand{\Ker}{\mathrm{Ker}\,}
\newcommand{\Hom}{\mathrm{Hom}\,}
\newcommand{\ttM}{\widetilde{\widetilde{M}}}
\newcommand{\ttN}{\widetilde{\widetilde{N}}}
\newcommand{\id}{\mathrm{id}}
\author{MASUDA Toshihiko \\
 Department of Mathematics, Kochi University, \\
2-5-1 Akebono-cho, Kochi, 780-8520, JAPAN}
\title{A note on group actions on subfactors}
\date{}
\begin{document}
\maketitle
\begin{abstract}
 We construct approximately inner actions of
 discrete amenable groups on strongly amenable subfactors of type
 II$_1$ with given invariants, and obtain classification results 
 under some conditions. 
We also study  the lifting of the relative $\chi$ group.
\end{abstract}

\section{Introduction}
In the theory of subfactors initiated by V. F. R. Jones in
\cite{J-ind}, the analysis of  automorphisms and group actions on
subfactors has been done by many people. 
We refer to \cite{CK},
\cite{EK-Hecke}, \cite{Go}, \cite{Go2},  \cite{Go3}, \cite{Kw-cnt},
\cite{Kw-orb}, \cite{Kw-kappa},
\cite{Kw-appr}, \cite{Ko1}, \cite{Ko-MSJ}, \cite{Li-auto},
\cite{Li-comm}, \cite{Po-act}, \cite{Xu-Orbi}. Also see \cite[Chapter
15]{EK-book}. 

In \cite{M-class}, 
we classified approximately inner actions of discrete
amenable groups on strongly amenable subfactors of type II$_1$ 
by the characteristic invariant and the $\nu$ invariant under 
some assumptions. Among these assumptions, the most important one is the
triviality of the algebraic $\kappa$ invariant. When the
algebraic $\kappa$ invariant  is trivial, 
we can classify approximately inner actions
completely. Hence we have to investigate the case when the algebraic
$\kappa$ invariant is not trivial. 
In this case, 
we do not know whether there
exist actions with given invariants or not. 
Hence what we should do first is to find a systematic
way to construct actions with given invariants. 
Note that if the algebraic $\kappa$ invariant is trivial, our
characteristic invariant is exactly same as original one in
\cite{J-act}, but if the algebraic $\kappa$ invariant is not trivial,  
our characteristic invariant may be different from the usual one, and
this makes classification more difficult.

In this paper, we construct actions of discrete amenable groups on
strongly amenable subfactors of type II$_1$ with given
invariants, and classify actions under an extra assumption on the $\nu$
invariant. 
(We emphasis that we never assume the triviality of the algebraic
$\kappa$ invariant.) 
The most essential assumption in our theory is that 
extensibility of 
the $\nu$ invariant to a homomorphism from a whole group.   
This assumption is similar to that of \cite[Theorem 20]{KwST}. 
In \cite{KwST}, 
Kawahigashi, Sutherland and Takesaki have 
classified the actions of a discrete abelian group $G$ on the injective type
III$_1$ factor. The modular invariant $\nu$ appears as the cocycle conjugacy
invariant, and this is a homomorphism from a subgroup of $G$ to $\mathbf{R}$.
Essential fact in their proof is that $\nu$  can be extended to
a homomorphism of $G$ due to the divisibility of $\mathbf{R}$. 
(Originally this idea was due to Connes. See \cite[pp.466]{Co-suv}.)

In subfactor case, we can not expect such property for the $\nu$
invariant generally. But 
if we assume the extensibility of 
the $\nu$ invariant, our proof goes well as in the proof of
 \cite[Theorem 20]{KwST}. 
We remark that our results can be viewed as the generalization of
\cite[Theorem 4.1]{Kw-appr}. 

In appendix, we discuss the lifting of $\chi_a(M,N)$ since
 we fix one lifting of $\chi_a(M,N)$ to
define characteristic invariant. 
\smallskip\\
\textbf{Acknowledgements.} The author is grateful to Professor M. Izumi for 
informing about \cite[Proposition 3.2]{Ol-P-Tak}, Professor Y. Kawahigashi for
various comments on this paper. The author is supported by Grant-in-Aid
for Scientific Research, Ministry of Education (Japan).

\section{Main results}
We use notations in \cite{M-class} freely.  

First we recall the definition of cocycle conjugacy invariants for
actions considered in \cite{M-class}. 
Let $N \subset M$ be a subfactor of type II$_1$ with finite index, 
$G$ a discrete group, 
and $\alpha$ an action of $G$ on $N\subset M$. 
Throughout this paper, we always make the following assumptions on
$N\subset M$. 
\smallskip\\
(A1) $N\subset M$ is extremal, \\
(A2) $N\subset M$ and $M\subset M_1$ have the trivial normalizer,\\
(A3) $\Ker\Phi=\Aut(M,N)$,\\
(A4) $\chi_a(M,N)$ is a finite group,\\
(A5) there exists a lifting $\sigma$ of $\chi_a(M,N)$ to
$\Aut(M,N)$. 

By (A3), every action has trivial Loi invariant. 
Note that we have many classes of subfactors satisfying the above
assumptions, e.g, Jones subfactor with principal graph
$A_{2n+1}$ in \cite{J-ind},
or subfactors coming from Hecke algebras in \cite{Wen-Hecke}. 
(Also see \cite{Kw-orb} and \cite{EK-Hecke}. )

By \cite[Theorem 3.1]{M-class}, we have a Connes-Radon-Nikodym type
cocycle $u_{\alpha,\sigma}\in U(N)$ for every $\alpha\in\Ker\Phi$ and $\sigma\in
\Cntr(M,N)$. The algebraic $\kappa$ invariant $\kappa_a$ is defined by 
$\kappa_a(h,k)=u_{\sigma_k,\sigma_h}^*$ for $h,k\in\chi_a(M,N)$.
We can easily verify that 
$\kappa_a$ is a bicharacter of 
$\chi_a(M,N)$.

For an action $\alpha$ of $G$, we get cocycle conjugacy invariants
in the following way. The first invariant is a normal subgroup
$H_\alpha \subset G$, which is a non-strongly outer part of $\alpha$. 
Then we get a $G$-equivariant homomorphism 
$\nu_\alpha$ from $H_\alpha$ to $\chi_a(M,N)$ by 
$\nu_\alpha(h)=[\alpha_h]$. This $\nu_\alpha$ is the second cocycle
conjugacy invariant, and we call this the $\nu$ invariant.
By (A5) $\alpha_h$ has the form
$\alpha_h=\Ad v_h \sigma_{\nu_\alpha(h)}$ for some unitary $v_h\in
U(N)$. Then we get a 
characteristic invariant $\Lambda(\alpha)=[\lambda_\alpha,\mu_\alpha]
\in\Lambda(G,H_\alpha|\kappa_a)$ by the following equations for $g\in
G$, $h,k\in H_\alpha$.
$$\alpha_g(v_{g^{-1}hg})u_{\alpha_g,\sigma_{\nu_\alpha(h)}}
=\lambda_\alpha(g,h)v_h,
\quad v_h\sigma_{\nu_\alpha(h)}(v_k)=\mu_\alpha(h,k)v_{hk}.$$   
The pair $\lambda(g,h)$ and $\mu(h,k)$ satisfy the following relations
for $h,k,l\in H$ and $g,g_1,g_2 \in G$. \\
(1) $\mu(h,k)\mu(hk,l)=\mu(k,l)\mu(h,kl),$  \\
(2) $\lambda(g_1g_2,h)=\lambda(g_1,h)\lambda(g_2,g_1^{-1}hg_1),$  \\
(3) $\lambda(g,hk)\overline{\lambda(g,h)\lambda(g,l)}=
\mu(h,k)\overline{\mu(g^{-1}hg,g^{-1}kg)}$ \\
(4) $\lambda(h,k)=
\mu(h,h^{-1}kh)\overline{\mu(k,h)\kappa_a(\nu(k),\nu(h))},$ \\
(5) $\lambda(e,h)=\lambda(g,e)=\mu(e,k)=\mu(h,e)=1.$ \\
 
Equation (4) shows the difference between the usual characteristic
invariant and our characteristic invariant. This definition of $\lambda$
and $\mu$  depends on the choice of $v_h$.  To get rid of this
dependence we have to define suitable equivalence relation for
$(\lambda,\mu)$. On this point see \cite{M-class}. 

Conversely for a given normal subgroup $H \subset G$, 
$[\lambda,\mu]\in\Lambda(G,H|\kappa_a)$ and
$\nu\in \Hom_G(H,\chi_a(M,N))$, we will construct an action $\alpha$ with 
$H_\alpha=H$, $\Lambda(\alpha)=[\lambda,\mu]$ and $\nu_\alpha=\nu$ in
the following proposition.

\begin{prop}\label{prop:model}
 Let $N\subset M$ be a strongly amenable subfactor of type II$_1$,  
 $G$ a discrete amenable group.
 Assume that $\nu$ can be extended to the homomorphism from $G$. 
 Then for every $[\lambda,\mu]\in\Lambda(G,H|\kappa_a)$ and $\nu$,
 there exists an action $\alpha$ of $G$ with $H_\alpha=H$,
 $\Lambda(\alpha)=[\lambda,\mu]$ and $\nu_\alpha=\nu$.
\end{prop}
\textbf{Proof.} By assumptions, we have an extension of $\nu$ from $G$
to $\chi_a(M,N)$, which we denote $\nu$ again.
Hence $g\to\sigma_{\nu(g)}$ is an action of $G$ on
$N\subset M$. Let $\kappa_a$ be the algebraic  $\kappa$ invariant 
for $N\subset M$, 
and set $\lambda'(g,n):=\kappa_a(\nu(n),\nu(g))\lambda(g,n)$. 
Then it is easy to verify that  
$[\lambda',\mu]$ is in $\Lambda(G,H)$, that is, $[\lambda',\mu]$ is a 
usual characteristic
invariant. Let $m$ be an action of $G$ on the injective type II$_1$
factor $R_0$ with the characteristic
invariant $[\lambda',\mu]$. Define an action $\alpha$ of $G$ by 
$\alpha_g:=\sigma_{\nu(g)}\otimes m_g$. 
Then this $\alpha$ is a desired one. \hfill$\Box$  

On classification of actions, we have the following result.
\begin{thm}\label{thm:class}
Let $N\subset M$, $G$ be as  in the previous proposition. 
Let $\alpha$ and $\beta$ 
be approximately inner actions of $G$. Assume $\nu_\alpha$ can be
extended to a homomorphism from G. 
Then $\alpha$ and $\beta$ are stably conjugate if 
$H_\alpha=H_\beta$,    
$\Lambda(\alpha)=\Lambda(\beta)$ and $\nu_\alpha=\nu_\beta$ hold. 
\end{thm}
\textbf{Proof.}
Set $K:=\chi_a(M,N)$. 
Let $\widetilde{\alpha}$ be an extension of $\alpha$ on $N\rtimes_\sigma
K\subset  M\rtimes_\sigma K$ defined in \cite{M-class}, 
and $\widetilde{\widetilde{\alpha}}$ be a
natural extension of $\widetilde{\alpha}$ on $
\ttN \subset \ttM :=N\rtimes_\sigma
K\rtimes_{\widehat{\sigma}}\widehat{K}\subset M \rtimes_\sigma K
\rtimes_{\widehat{\sigma}}\widehat{K}$. 
Let $w_k$ be an implementing unitary of $\sigma$ in $M\rtimes_\sigma K$,
and $v_p$ be an implementing unitary of $\widehat{\sigma}$ in $M\rtimes_\sigma
K\rtimes_{\widehat{\sigma}}\widehat{K}$. 
Then by the definition of $\widetilde{\widetilde{\alpha}}$, we have 
$\widetilde{\widetilde{\alpha}}_g(x)=\alpha_g(x)$, 
$\widetilde{\widetilde{\alpha}}_g(w_k)=u_{\alpha_g,\sigma_k}w_k$ and 
$\widetilde{\widetilde{\alpha}}_g(v_p)=v_p$ for $x\in M$, $k\in K$ and
$p\in\widehat{K}$.  
On the other hand, the second dual action $\widehat{\widehat{\sigma}}$
of $\sigma$ satisfies $\widehat{\widehat{\sigma}}_k=\mathrm{id} $
on $M\rtimes_\sigma K$ and $\widehat{\widehat{\sigma}}_k(v_p)=
\overline{\langle k,p\rangle} v_p$ for $p\in \widehat{K}$. 

Takesaki duality theorem says that  
$\ttN\subset \ttM$ is isomorphic to 
$N\otimes B(l^2(K)) \subset M \otimes B(l^2(K))$ 
via an isomorphism $\Psi$ satisfying the following. 
\begin{eqnarray*}
(1)&& (\Psi(\pi_{\widehat{\sigma}}\circ\pi_\sigma(a))\xi)(k)
    =\sigma_k^{-1}(a)\xi(k), \\
(2)&& (\Psi(\pi_{\widehat{\sigma}}(w_l))\xi)(k)=\xi(l^{-1}k), \\
(3)&& (\Psi(v_p)\xi)(k)=\overline{\langle k,p\rangle}\xi(k),
\end{eqnarray*}
where $\pi_\sigma$ is an embedding of $M$ into $M\rtimes_\sigma K$, and 
$\pi_{\widehat{\sigma}}$ is an embedding of $M\rtimes_\sigma K$ into 
$M\rtimes_\sigma K \rtimes_{\widehat{\sigma}} \widehat{K}$. 

Define a unitary $c_g\in N\otimes B(l^2(K))$ by
$(c_g\xi)(k):=u_{\alpha_g,\sigma_k^{-1}}^*\xi(k)$. Since $c_g$ commutes
with elements in $N'\otimes \mathbf{C}1$, $c_g$ is indeed in $N\otimes
B(l^2(K))$. Moreover since we have    
\begin{eqnarray*}
 (c_g\alpha_g\otimes id(c_h)\xi)(k) &=& u_{\alpha_g,\sigma_k^{-1}}^* 
         \alpha_g(u_{\alpha_h,\sigma^{-1}_k}^*)\xi(k) \\
 &=&u_{\alpha_{gh},\sigma_k^{-1}}^*\xi(k) \\
 &=& (c_{gh}\xi)(k),
\end{eqnarray*}
$c_g$ is an $\alpha\otimes \id $ cocycle. 
Then as in  the argument in \cite[Section 5]{Ko-MSJ}, 
it is shown that $\Psi\circ\widetilde{\widetilde{\alpha}}_g\circ \Psi^{-1}= 
\Ad c_g (\alpha_g\otimes \id)$ holds.  

On the other hand we have $\Psi\circ\widehat{\widehat{\sigma}}_k\circ
\Psi^{-1}=\sigma_k\otimes \Ad \rho_k^{-1}$, where $\rho$ is a left
regular representation of $K$.

Here we consider the Connes-Radon-Nikodym type cocycle for 
$\Ad c_g \alpha_g\otimes \id$ and $\sigma_k\otimes\Ad \rho_k^{-1}$.
Take $0\ne a\in M_n$ with $\sigma_k(x)a=ax$ for every $x\in M$. 
By \cite[Theorem 3.1]{M-class}, $\alpha_g(a)=u_{\alpha_g,\sigma_k}a$
holds. It is obvious that 
$\sigma_k\otimes \Ad \rho^{-1}_k (x)(a\otimes \rho^{-1}_k)=(a\otimes
\rho^{-1}_k) x$ holds for every $M\otimes B(l^2(K))$.
Here we have the following. 
\begin{eqnarray*}
 (\Ad c_g (\alpha_g\otimes \id)(a\otimes \rho^{-1}_k)\xi)(l)&=&
 (c_g(\alpha_g(a)\otimes \rho_k^{-1})c_g^*\xi)(l)\\
 &=&u_{\alpha_g,\sigma_l^{-1}}^*\alpha_g(a)(c_g^*\xi)(kl)\\
 &=&u_{\alpha_g,\sigma_l^{-1}}^*u_{\alpha_g,\sigma_k}a
  u_{\alpha_g,\sigma_{kl}^{-1}}\xi(kl)\\
 &=&u_{\alpha_g,\sigma_l^{-1}}^*u_{\alpha_g,\sigma_k}\sigma_k(u_{\alpha_g,\sigma_{kl}^{-1}})a\xi(kl) \\
 &=&u_{\alpha_g,\sigma_l^{-1}}^*u_{\alpha_g,\sigma_l^{-1}}a\xi(kl)\\
 &=&(a\otimes \rho_k^{-1}\xi)(l).
\end{eqnarray*}

By \cite[Theorem 3.1]{M-class}, the above equality implies $u_{\Ad
c_g(\alpha\otimes \id), \sigma_k\otimes \Ad \rho_k^{-1}}=1$ holds for
every $g\in G$ and $k\in K$. 
Hence by replacing $\alpha$ and $\beta$ if necessary, we may assume that 
$u_{\alpha_g,\sigma_k}=1$ and 
$u_{\beta_g,\sigma_k}=1$ hold for every $g\in G$ and $k\in K$. 
This especially implies $\alpha_g\sigma_k=\sigma_k\alpha_g$ and 
$\beta_g\sigma_k=\sigma_k\beta_g$.

Define two new actions ${\bar \alpha}$ and ${\bar \beta}$ by
${\bar \alpha}_g:=\alpha_g\sigma_{\nu(g)}^{-1}$ and 
${\bar \beta}_g:=\beta_g\sigma_{\nu(g)}^{-1}$. 
Since $\alpha$ and 
$\beta$ commute with $\sigma$, 
${\bar \alpha}$ and 
${\bar \beta}$ are indeed actions of $G$. 
By construction of ${\bar \alpha}$ and ${\bar \beta}$, it is easy to see
$H_\alpha={\bar \alpha}^{-1}(\Cnt(M,N))={\bar
\alpha}^{-1}(\Int(M,N))={\bar \beta}^{-1}(\Cnt(M,N))={\bar
\beta}^{-1}(\Int(M,N))$. 

Next we compute $\Lambda({\bar \alpha})$. 
Take $m\in H$ and $v_m\in U(N)$ with $\alpha_m=\Ad v_m\sigma_{\nu(m)}$. 
In this case, we have $\Ad {\bar\alpha}_m=\Ad v_m$ for $m\in H$.
Moreover since 
\begin{eqnarray*}
 1&=&u_{\alpha_h,\sigma_k}\\
  &=&u_{\Ad v_m\sigma_{\nu(m)},\sigma_k } \\
  &=&\Ad v_m (u_{\sigma_{\nu(m)},\sigma_k}) u_{\Ad v_m,\sigma_k}\\
  &=&\overline{\kappa_a(k,\nu(m))}v_m\sigma_k(v_m^*)
\end{eqnarray*}
holds, we have $\sigma_k(v_m)=\overline{\kappa_a(k,\nu(m))}v_m$.

First we compute $\lambda_{\bar{\alpha}}$.  
Since we have $u_{\alpha_g,\sigma_k}=1$,
$\alpha_g(v_{g^{-1}ng})=\lambda_\alpha(g,n)v_n$ holds by the definition
of $\lambda_\alpha$. 
Then we get
\begin{eqnarray*}
 \bar{\alpha}_g(v_{g^{-1}ng})&=&\alpha_g\sigma_{\nu(g)}^{-1}(v_{g^{-1}ng})\\ 
 &=&\overline{\kappa_a(\nu(g)^{-1},\nu(n))}\alpha_g(v_{g^{-1}ng}) \\
 &=&\kappa_a(\nu(g),\nu(n))\lambda_\alpha(g,n)v_n,
\end{eqnarray*} 
and
$\lambda_{\bar\alpha}(g,n)=\kappa_a(\nu(g),\nu(n))\lambda_\alpha(g,n)$
holds. 
Next we compute $\mu_{\bar{\alpha}}$.   
By the definition of $\mu_\alpha$, we have
$v_m\sigma_{\nu(m)}(v_n)=\mu_\alpha(m,n)v_{mn}$, $m,n\in H$. Hence we get
$v_mv_n=\kappa_a(\nu(m),\nu(n))\mu_\alpha(m,n)v_{mn}$ and
consequently 
$\mu_{\bar{\alpha}}(m,n)=\kappa_a(\nu(m),\nu(n))\mu_\alpha(m,n)$.

Similar computation holds for ${\bar \beta}$, and by the assumption 
$\Lambda(\alpha)=\Lambda(\beta)$, we get
$\Lambda(\bar{\alpha})=\Lambda(\bar{\beta})$. 
Hence ${\bar \alpha}$ and ${\bar \beta}$ are cocycle conjugate 
by \cite[Theorem 5.1]{M-class}. 
Then ${\bar \alpha}$ and ${\bar \beta}$ are stably conjugate, hence 
there exists an automorphism $\theta\in \Aut(M\otimes B(l^2(G)),
N\otimes B(l^2(G)))$ with $\theta\circ({\bar \beta}_g\otimes \Ad
\varrho_g) \circ \theta^{-1}={\bar \alpha}_g$, where $\varrho$ is a
right regular representation of $G$. 

To prove the main theorem, we need the following proposition. 
In the following proposition, $N\subset M$ can be an arbitrary subfactor
of finite index.
\begin{prop}\label{prop:fix}
Let $\sigma$ be a non-strongly outer automorphism. 
Take $0\ne a\in M_n $ such that $\sigma(x)a=ax$ holds for every $x\in M$. 
Then $v\in M$ is in $M^\sigma$ if and only if $va=av$ holds.
\end{prop}
\textbf{Proof.} First assume that $v\in M^\sigma$. Then we have 
$va=\sigma(v)a=av$. Conversely assume that $va=av$ holds. Then
$\sigma(v)a=av=va$ holds. Hence we have $\sigma(v)aa^*=vaa^*$. Here
$aa^* $ is in $M'\cap M_n$. Let $E$ be the minimal conditional
expectation from $M_n$ onto $M$. 
Then we get $\sigma(v)E(aa^*)=E(\sigma(v)aa^*)=E(vaa^*)=vE(aa^*)$,
and $E(aa^*)\in M\cap M'=\mathbf{C}$. Since $a$ is not zero, $E(aa^*)$
is a non-zero scalar. Hence $v$ is in $M^\sigma$. \hfill$\Box$ 
\smallskip \\  
\textbf{Remark.} The above proposition can be regarded as a 
subfactor-analogue of the characterization of the centralizer of type III
factors. Namely let $M$ be a type III factor, $\phi$ a faithful normal
state of $M$. Then $a$ is in $M_\phi$ if and only if $[\phi, a]=0$.
\smallskip

We continue the proof of Theorem \ref{thm:class}. 
Since an outer action of a finite group is stable, 
we can find a unitary $w\in N\otimes B(l^2(G))$ such that 
$w^*\sigma_k\otimes \id(w)=u_{\theta,\sigma_k\otimes\id}$. 
Hence $\theta\circ\sigma_k\otimes\id\circ
\theta^{-1}= \Ad u_{\theta, \sigma_k\otimes \id}\circ \sigma_k\otimes \id=
\Ad w^*\circ \sigma_k\otimes \id\circ \Ad w$ holds. If we can prove
that $w{\bar \alpha}_g\otimes \Ad \varrho_g(w^*)$ is in 
$(M\otimes B(l^2(G)))^K$, 
then $w{\bar \alpha}_g\otimes \Ad \varrho_g(w^*)$ is an 
$\alpha\otimes\Ad\varrho={\bar \alpha}\sigma\otimes \Ad \varrho$ 
cocycle and 
\begin{eqnarray*}
{\bar \alpha}_g\otimes\Ad\varrho_g \circ
\theta\circ\sigma_{\nu(g)}\otimes \id\circ\theta^{-1}&=&
 {\bar \alpha}_g\otimes \Ad\varrho_g\circ \Ad w^*\circ\sigma_{\nu(g)}
 \otimes \id \circ\Ad w \\
&=& \Ad ({\bar \alpha}\otimes 
\Ad\varrho_g(w^*)){\bar \alpha}_g\sigma_{\nu(g)}\otimes\Ad
  \varrho_g \circ \Ad w \\
&=& \Ad w^* \circ \Ad (w{\bar \alpha_g}\otimes\Ad\varrho_g(w^*))
  \alpha_g \otimes \Ad\varrho_g\circ \Ad w 
\end{eqnarray*}
holds and we have the following. 
\begin{eqnarray*}
 \alpha\otimes\Ad\varrho&=&{\bar \alpha}\otimes\Ad\varrho\circ
      \sigma\otimes\id \\  
      &\sim& 
      {\bar \alpha}\otimes\Ad\varrho\circ\theta\circ\sigma\otimes\id
      \circ\theta^{-1} \\ 
       &=&\theta\circ \bar{\beta}\otimes\Ad\varrho\circ \theta^{-1}\theta 
       \circ\sigma\otimes\id\circ\theta^{-1} \\
       &=&\theta\circ \bar{\beta}\sigma\otimes\Ad\varrho\circ \theta\\
       &=&\theta\circ \beta\otimes\Ad\varrho \circ\theta^{-1}. 
\end{eqnarray*}
Hence $\alpha$ and $\beta$ are stably conjugate. 
So we only have to prove that $w{\bar \alpha}_g\otimes\Ad\varrho_g(w^*) $
is in $(M\otimes B(l^2(G)))^K$.

It is easy to see 
$u_{\alpha_g\otimes\Ad\varrho_g,\sigma_k\otimes\id}=1$, hence 
$u_{\bar{\alpha}_g\otimes\Ad\varrho_g,\sigma_k\otimes\id}=
 u_{\alpha_g\sigma_{\nu(g)}^{-1}\otimes\Ad\varrho_g,
 \sigma_k\otimes\id}=\kappa_a(k,\nu(g))$ holds.
In the same way, we can see 
$u_{\bar{\beta}_g\otimes\Ad\varrho_g,\sigma_k\otimes\id}=\kappa_a(k,\nu(g))$. 
Hence 
\begin{eqnarray*}
 u_{\bar{\alpha}_g\otimes\Ad\varrho_g, 
  \Ad w^*\circ\sigma_k\otimes \id\circ \Ad w}&=&
  u_{\theta\circ\bar{\beta}_g\otimes\Ad\varrho_g\circ\theta^{-1}, 
  \theta\circ\sigma_k\otimes \id\circ \theta^{-1}}  \\
 &=&\theta(u_{\bar{\beta}_g\otimes\Ad\varrho_g, 
  \sigma_k\otimes \id})\\
 &=&\kappa_a(k,\nu(g)) 
\end{eqnarray*}
holds.
Take $0\ne a\in M_n\otimes B(l^2(G))$ such that 
$\sigma_k\otimes\id(x)a=ax$ holds for every $x\in M\otimes B(l^2(G))$.
Then ${\bar\alpha}_g\otimes\Ad \varrho_g(a)=\kappa_a(k,\nu(g))a$ holds
by \cite[Theorem 3.1]{M-class}. 
Since $\Ad w^*\circ\sigma_k\otimes\id\circ\Ad w(x)
w^* aw=  w^*aw x$, we also have ${\bar
\alpha}_g\otimes\Ad\varrho_g(w^*aw)=\kappa_a(k,\nu(g))w^*aw$. 
From these two equalities, we get 
$\bar{\alpha}_g\otimes\Ad\varrho_g(w^*)a\bar{\alpha}_g\otimes \Ad
\varrho_g(w)=w^*aw$. Hence $w\bar{\alpha}_g\otimes\Ad\varrho_g(w^*)$ satisfies
the condition in Proposition \ref{prop:fix}, 
$w\bar{\alpha}_g\otimes\Ad\varrho_g(w^*)$ is
in a fixed point algebra $(M\otimes B(l^2(G)))^K$. \hfill$\Box$

\begin{cor}
If $G$ is a finite group in Theorem \ref{thm:class}, 
then $\alpha$ and $\beta$ are cocycle conjugate if and only if
 $H_\alpha=H_\beta$, $\Lambda(\alpha)=\Lambda(\beta)$ and
 $\nu_\alpha=\nu_\beta$ hold.  
\end{cor}
\textbf{Proof.} By Theorem \ref{thm:class}, $\alpha\otimes\Ad\varrho$ and
$\beta\otimes \Ad\varrho$ are conjugate.
But in the same way as in the
proof of \cite[Lemma 6.5]{KtST}, we can prove that $\alpha$ is cocycle
conjugate to $\alpha\otimes m$, where $m$ is an outer action of $G/H$ on
the injective type II$_1$ factor $R_0$ and we regard $m$ as an action of
$G$ in the natural way.
Hence $\alpha$ and $\beta$ are cocycle conjugate since 
$m$ and $m\otimes\Ad\varrho$ are cocycle conjugate. \hfill$\Box$

In the rest of this paper, we treat examples which satisfy the
assumption in Theorem \ref{thm:class}. The first example is 
taken from \cite[Theorem 4.1]{Kw-appr}.

\begin{exam}
 \upshape
 We consider the case $G=\mathbf{Z}$. Take $\alpha\in\Aut(M,N)$. 
 Let $p$ be a strongly outer period of $\alpha$. 
  Set $\sigma:=\alpha^{p}$. 
 Then $\nu_\alpha$ is given by $\nu_\alpha(pm)=[\sigma^m]$.
 Let $n$ be an outer period of $\sigma$.  
 Here assume $(p,n)=1$. Then we can find $k,l\in\mathbf{Z}$ such
 that $pk+nl=1$. Set $\nu(g):=[\sigma^{gk}]$. Then we have
 $\nu(p)=[\sigma^{pk}]=[\sigma^{-nl+1}]=[\sigma]$, hence $\nu_\alpha$
 can be extended to a homomorphism from $\mathbf{Z}$.
\end{exam}

\begin{exam}
 \upshape
 Assume that $G$ is of the form $G=H_\alpha\rtimes K$. 
  For $(h,k)\in G=H_\alpha\rtimes K$, 
 define $\nu(h,k):=\nu_\alpha(h)$. Then by using the fact 
 $\nu_\alpha(knk^{-1})=\nu_\alpha(n)$, we get 
 \begin{eqnarray*}
  \nu((h_1,k_1)(h_2,k_2))&=&\nu(h_1k_1h_2k_1^{-1}, k_1k_2) \\
  &=&\nu_\alpha(h_1k_1h_2k_1^{-1}) \\
  &=&\nu_\alpha(h_1)\nu_\alpha(k_1h_2k_1^{-1}) \\
  &=&\nu(h_1,k_1)\nu(h_2,k_2).
 \end{eqnarray*}
Hence we can extend $\nu_\alpha$ to the homomorphism from $G$, 
and we can apply main theorem. 

\end{exam}

\appendix
\section{On liftings of the relative $\chi$ group}
In \cite{M-class} and this paper, we fixed a lifting of $\chi_a(M,N)$
to $\Aut(M,N)$ for classification of group actions on subfactors.  
But it seems that there are no way to choose the natural lifting. In this
appendix, we show that there exists a natural choice of a lifting by using the
algebraic $\kappa$ invariant. 

Take a lifting $\sigma$. Then the algebraic $\kappa$
invariant $\kappa_a(h,k)$ is defined as 
$\kappa_a(h,k):=u_{\sigma_k,\sigma_h}^*$. To specify
$\sigma$, we denote this $\kappa_a$ by
$\kappa_a^\sigma(h,k)$. 

Fix $\sigma$ and take another lifting
$\tilde{\sigma}$. Then we can find a unitary $u_h\in U(N)$ with $\Ad
u_h\sigma_h=\tilde{\sigma}_h$. Since $\tilde{\sigma}$ is a lifting,
there exists a 2-cocycle $\mu(h,k)\in Z^2(\chi_a(M,N),\mathbf{T})$ with
$u_h\sigma_h(u_k)=\mu(h,k)u_{hk}$. Here we compute
$\kappa_a^{\tilde{\sigma}}$.  Then  
\begin{eqnarray*}
\overline{\kappa_a^{\tilde{\sigma}}(h,k)}
&=&u_{\tilde{\sigma}_k,\tilde{\sigma}_h} \\
&=&u_{\Ad u_k\sigma_k,\Ad u_h\sigma_h} \\
&=&u_{\Ad u_k\sigma_k, \Ad u_h}\Ad u_h(u_{\Ad u_k\sigma_k, \sigma_h}) \\ 
&=& u_k\sigma_k(u_h)u_k^*u_h^* 
  \Ad u_h(\Ad u_k(u_{\sigma_k,\sigma_h})u_{\Ad u_k,\sigma_h}) \\
&=& u_k\sigma_k(u_h)u_k^*u_{\sigma_k,\sigma_h}u_k\sigma_h(u_k^*)u_h^* \\
&=& \overline{\kappa_a^\sigma(h,k)}\overline{\mu(h,k)}\mu(k,h)
\end{eqnarray*}
holds, hence we get
$\kappa_a^{\tilde{\sigma}}(h,k)=
\kappa_a^\sigma(h,k)\mu(h,k)\overline{\mu(k,h)}$.    

Here assume $\kappa_a^{\tilde{\sigma}}=\kappa_a^\sigma$. Then we
get $\mu(h,k)=\mu(k,h)$. By \cite[Proposition 3.2 ]{Ol-P-Tak}, $\mu$ is
a coboundary, so we can choose $u_h$ as a $\sigma$-cocycle. Hence
we get the following proposition. 

\begin{prop}\label{prop:lift}
 Let $\sigma$ and $\tilde{\sigma}$ be  liftings of $\chi_a(M,N)$ to
 $\Aut(M,N)$. If $\kappa_a^\sigma=\kappa_a^{\tilde{\sigma}}$ holds, then 
 $\tilde{\sigma}$ is a cocycle perturbation of $\sigma$. 
\end{prop}

By Proposition \ref{prop:lift}, we can find a lifting $\sigma$ up to cocycle
perturbation once we fix the algebraic $\kappa$ invariant. 

In the next proposition, 
we do not assume $\Ker\Phi=\Aut(M,N )$. 
Every $\theta \in \Aut(M,N)$ induces an automorphism
$\chi_a(\theta)$ of
$\chi_a(M,N)$ by
$\chi_a(\theta)([\sigma]):=[\theta\circ\sigma\circ\theta^{-1}]$.   

\begin{prop}
Let $\sigma$ be a lifting of $\chi_a(M,N)$ to $\Aut(M,N)$. 
Assume that
 $\kappa_a^\sigma(h,k)=\kappa_a^\sigma(\chi_a(\theta)(h),\chi_a(\theta)(k))$
 holds for every $\theta\in \Aut (M,N)$.  
Then there exists a $\sigma_{\chi_a(\theta)(\cdot)}$-cocycle $w_h$ such that 
$\theta\circ\sigma_h\circ\theta^{-1}=\Ad w_h 
\sigma_{\chi_a^\sigma(\theta)(h)}$ holds. 

\end{prop}
\textbf{Proof.} Take a unitary $w_h$ with $\Ad
w_h\sigma_{\chi_a(\theta)(h)}=\theta\circ\sigma_h\circ\theta^{-1}$. 
Then there exists a 2-cocycle $\mu(h,k)$ satisfying
$w_h\sigma_{\chi_a(\theta)(h)}(w_k)=\mu(h,k)w_{hk}$. 
On one hand, we have $u_{\theta\circ\sigma_h\circ\theta^{-1},
\theta\circ\sigma_k\circ\theta^{-1}}^*=
\theta(u_{\sigma_h,\sigma_k})^*=\kappa_a^\sigma(k,h)$.
On the other hand, we have
\begin{eqnarray*}
u_{\theta\circ\sigma_h\circ \theta^{-1},
\theta\circ\sigma_k\circ\theta^{-1}}^*
&=&u_{\Ad
w_h\sigma_{\chi_a(\theta)(h)},\Ad w_k\sigma_{\chi_a(\theta)(k)}}^* \\
&=&\kappa_a^\sigma(\chi_a(\theta)(k),\chi_a(\theta)(h))
\mu(k,h)\overline{\mu(h,k)}. 
\end{eqnarray*}
By assumption on $\kappa_a^\sigma$, we can choose $w_h$ as a cocycle
as the same reason in the proof of Proposition \ref{prop:lift}.  
\hfill$\Box$ 

The assumption on $\kappa_a$ in the above proposition is satisfied when
(1) $\kappa_a$ is trivial, (2) $\chi_a(M,N)$ is a cyclic group. 
The former is obvious. The reason of latter is following. It is easy to
see
$\kappa_a^\sigma(h,h)=\kappa_a^\sigma(\chi_a(\theta)(h),\chi_a(\theta)(h))$
holds from the above computation. 
If $\chi_a(M,N)$ is cyclic and $g$ is a generator of
$\chi_a(M,N)$, then
\begin{eqnarray*}
 \kappa_a^\sigma(g^m,g^n)&=&\kappa_a^\sigma(g,g)^{mn}\\
 &=&\kappa_a^\sigma(\chi_a(\theta)(g),\chi_a(\theta)(g))^{mn}\\
 &=&\kappa_a^\sigma(\chi_a(\theta)(g^m),\chi_a(\theta)(g^n))
\end{eqnarray*}
holds.

\end{document}